\documentclass[10pt]{article}
\newtheorem{thm}{Theorem}
\newtheorem{lema}[thm]{Lemma}
\newtheorem{propo}[thm]{Proposition}
\newtheorem{coro}[thm]{Corollary}

\newcommand{\Z}{\hbox{\bf Z}}

\newcommand{\beeq}{\begin{eqnarray*}}
\newcommand{\eneq}{\end{eqnarray*}}
\newcommand{\proof}{\noindent {\it Proof.\hspace{4mm}}}
\newcommand{\qfd}{\hfill $\fbox{}$\vspace{4mm}}

\newcommand\ZZ{{{\rm Z}\kern-.28em{\rm Z}}}
\def\newpic#1{%
\def\emline##1##2##3##4##5##6{%
         \put(##1,##2){\special{em:point #1##3}}%
         \put(##4,##5){\special{em:point #1##6}}%
         \special{em:line #1##3,#1##6}}}
\newpic{}
\def\emline#1#2#3#4#5#6{%
          \put(#1,#2){\special{em:moveto}}%
          \put(#4,#5){\special{em:lineto}}}
\def\newpic#1{}
\title{SQS-graphs of Solov'eva-Phelps codes}
\author{Italo J. Dejter
\\ University of Puerto Rico \\ Rio Piedras, PR 00931-3355 \\ ijdejter@uprrp.edu
}

\date{}
\begin{document}
\maketitle
\centerline{\it Dedicated to Charles C. Lindner's in his 70th anniversary}

\begin{abstract}
A binary extended 1-perfect code $\mathcal C$ folds over its kernel
via the Steiner quadruple systems associated with its codewords. The
resulting folding, proposed as a graph invariant for $\mathcal C$,
distinguishes among the 361 nonlinear codes $\mathcal C$ of kernel
dimension $\kappa$ obtained via
Solov'eva-Phelps doubling construction, where $9\geq\kappa\geq 5$. Each of the 361 resulting
graphs has most of its nonloop edges expressible in terms of
lexicographically ordered quarters of products of classes from extended
1-perfect partitions of length 8 (as classified by Phelps)
and loops mostly expressible in terms of the
lines of the Fano plane.
\end{abstract}

\section{Preliminaries, objectives and plan}

We consider the $n$-cube $Q_n$ as the graph with
vertex set $F_2^n=\{0,1\}^n$ in which each two vertices
that differ in exactly one coordinate are joined by an edge. A perfect 1-error-correcting
code, or 1-perfect code, $C=C^r$ of length $n=2^r-1$, where
$0<r\in\Z$, is an independent vertex set of $Q_n$ such that each
vertex of $Q_n\setminus C$ is neighbor of exactly one vertex of $C$.
It follows that $C$ has distance 3 and $2^{n-r}$ vertices.

Each 1-perfect code $C=C^r$ of length $n=2^r-1$ can be extended by
adding an overall parity check. This yields an {\it extended} 1-{\it
perfect code} ${\mathcal C}={\mathcal C}^r$ of length $n+1=2^r$,
which is a subspace of even-weight words of $F_2^{n+1}$. The $n+1$
coordinates of the words of $F_2^{n+1}$ here are orderly indicated
$0,1,\ldots,n$.

For every $n=2^r-1$ such that $0<r\in\Z$ there is at least one
linear code $C^r$ as above and a corresponding linear extension, ${\mathcal C}^r$.
These codes are unique for every $r<4$. The situation changes for
$r\ge 4$. In fact, there are many nonlinear codes $C^4$ and
${\mathcal C}^4$ or length 15 and 16, respectively,
\cite{EV,Ph1,PhL,Rif,Sol1,Vas}.

The kernel $Ker(C)$ of a 1-perfect code $C$ of length $n$ is defined
as the largest subset $K\subseteq Q_n$ such that any vector in $K$
leaves $C$ invariant under translations, \cite{PhL}. In other words,
$x\in Q_n$ is in $Ker(C)$ if and only if $x+C=C$. If $C$  contains
the zero vector, then $Ker(C)\subseteq C$. In this case, $Ker(C)$ is
also the intersection of all maximal linear subcodes contained in
$C$. The kernel $Ker({\mathcal C})$ of an extended 1-perfect code
$\mathcal C$ of length $n+1$ is defined in a similar fashion in
$Q_{n+1}$.

A partition of $F_2^n$ into 1-perfect codes $C_0,C_1,\ldots,C_n$ is
said to be a 1-{\it perfect partition} $\{C_0,C_1,\ldots,C_n\}$ of
length $n$. The following result on {\it doubling construction}
of extended 1-per\-fect codes of length $2n+2$ is due to Solov'eva
\cite{Sol1} and Phelps \cite{Ph1}, so in this work they are called
{\it SP-codes}.

\begin{thm}{\rm\cite{Sol1,Ph1,Ph2}}
Given two extended 1-perfect partitions $\{C_0,C_1,$ $\ldots,C_n\}$
and $\{D_{n+1},D_{n+2},\ldots,D_{2n}\}$ of length $n+1$ and a
permutation $\sigma$ of $[0,n]=\{0,1,\ldots,n\}$, there exists a
1-perfect code $\mathcal C$ of length $2n+2$ given by ${\mathcal
C}=\bigcup_{i=0}^n\left\{(x,y)|x\in C_i,y\in
D_{n+1+\sigma(i)}\right\}.$
\end{thm}

Few invariants for 1-perfect codes $C$ have been proposed toward their
classification, namely: the rank of $C$ and the dimension of
$Ker(C)$ \cite{PhL}, the STS-graph $H(C)$ \cite{Dej1} and the
STS-graph $H_K(C)$ modulo $K=Ker(C)$ \cite{Dej2}.

In the present work, an invariant for extended 1-perfect codes
$\mathcal C$, referred to as the {\it SQS-graph} $H_K({\mathcal C})$
of $\mathcal C$ modulo $K=Ker({\mathcal C})$, is presented, and computed in particular for the SP-codes of length 16 and kernel dimension $\kappa$
such that $9\ge\kappa\ge 5$, allowing for us to distinguish successfully between each two of them
and to show in Theorem 5 that each such $H_K({\mathcal C})$ has its
{\it links} (nonloop edges) expressible in terms of products of
classes from extended 1-perfect
partitions $\{C_0,\ldots,C_7\}$ and
$\{D_8,\ldots,D_{15}\},$ and their loops mostly expressible in terms
of the lines of the Fano plane.

In \cite{Ph2}, Phelps found that there are exactly eleven 1-perfect
partitions of length 7, denoted $0,1,\ldots,10$. If two such
partitions are equivalent, then the corresponding extended
partitions are equivalent. However, the converse is false. In fact,
puncturing an extended 1-perfect partition at different coordinates
may result in nonequivalent 1-perfect partitions. Also, Phelps found
that there are just ten nonequivalent extended 1-perfect partitions
of length 8. In fact, partitions 2 and 7 in \cite{Ph2} have
equivalent extensions.

Phelps also found in \cite{Ph2} that there are exactly 963 extended
1-perfect codes of length 16 obtained via Theorem 1 applied to the cited partitions. The members in the
corresponding list of 963 SP-codes in \cite{Ph2} are referred below
according to their order of presentation, with numeric indications
(using three digits), from 001 (corresponding to the linear code) to
963,  (or from 1 to 963). This numeric indication is presented in an
additional final column in a copy of the listing of \cite{Ph2} that
can be retrieved from http://home.coqui.net/dejterij/963.txt. It
contains, for each one of its 963 lines: a reference number,
rank, dimension of the kernel, two numbers in
$\{0,1,2,3,4,5,6,8,9,10=a\}$ representing corresponding Phelps'
1-perfect source partition $\{C_0,\ldots,C_7\}$ and target partition
$\{D_8,\ldots,D_{15}\}$ and a permutation $\sigma$ as in Theorem 1.

In Section 2, $H_K({\mathcal C})$ is defined and subsequently
applied to the nonli\-near SP-codes of length 16 and kernel dimension
$\kappa\ge 5$, according to their classification in \cite{Ph2}. We
note that there are 361 nonlinear SP-codes of length 16 with
$\kappa\geq 5$, namely 2, 10, 18, 86 and 245 respectively with $\kappa=9,8,7,6$ and 5.
Section 3 accounts for the participating STS(15)-types
defined in Section 2 and in accordance with \cite{LeVan,MPR,WCC}. Section 4 accounts for those SP-codes that show
homogeneous behavior with respect to the involved Steiner
quadruple and triple systems, according to the results of Theorem 5, (Section 5), with which we
deal in the the rest of the paper.

\section{Definition of $H_K({\mathcal C})$}

A {\it Steiner quadruple system}, (or SQS), is an  ordered pair
$(V,B)$, where $V$ is a finite set and $B$ is a set of quadruples of
$V$ such that every triple of $V$ is a subset of exactly one
quadruple in $B$. A subset of $B$ will be said to be an {\it SQS-subset}.

The minimum-distance graph $M({\mathcal C})$ of an extended
1-perfect code $\mathcal C$ of length $n$ has $\mathcal C$ as its vertex set and
exactly one edge between each two vertices $v,w\in C$ whose Hamming
distance is $d(v,w)=4$. Each edge $vw$ of $M({\mathcal C})$ is
naturally labeled with the quadruple of coordinate indices
$i\in\{0,\ldots,n\}$ realizing $d(v,w)=4$. As a result, the labels
of the edges of $M({\mathcal C})$ incident to any particular vertex
$v$ constitute a Steiner quadruple system $S({\mathcal C},v)$ formed
by $n(n+1)(n-1)/24$ quadruples on the $n+1$ coordinate indices
which in our case, namely for $n+1=16$, totals 140
quadruples.

Given an edge $vw$ of $M({\mathcal C})$, its labeling quadruple is
denoted $s(vw)$. Each codeword $v$ of $\mathcal C$ is labeled by the
equivalence class ${\cal S}[v]$ of Steiner quadruple systems on $n$
elements corres\-ponding to $S({\mathcal C},v)$, called for short
the {\it SQS($n$)-type} ${\cal S}[v]$. We say that $M({\mathcal C})$
with all these vertex and edge labels is the {\it SQS-graph} of
$\mathcal C$. Let $L\subseteq K=Ker({\mathcal C})$ be a linear
subspace of $\mathcal C$. Clearly, $L$ partitions $\mathcal C$ into
classes $v+L$. ($w\in{\mathcal C}$ is in $v+L$ if and only if
$v-w\in L$). These classes $v+L$ are said to be the {\it classes} of
$\mathcal C$ mod $L$. The set they form can be taken as a quotient
set ${\mathcal C}/L$ of $\mathcal C$. The following three results
and accompanying comments are similar in nature to corresponding
results in \cite{Dej2}, but now the code $\mathcal C$ in their
statements is assumed to be an extended 1-perfect code.

\begin{lema}
Each $v+L\in{\mathcal C}/L$ can be assigned a well-defined Steiner
quadruple system $S({\mathcal C},v)$. \qfd
\end{lema}

Lemma 2 suggests the following `foldability' condition. If for any
two classes $u+L$ and $v+L$ of $\mathcal C$ mod $L$ with $d(u,v)=4$
realized by $s(uv)$ holds that for any $u'\in u+L$ there is a $v'\in
v+L$ with $d(u',v')=4$ and realized {\it exactly} by
$s(u'v')=s(uv)$, then we say that $\mathcal C$ is {\it foldable
over} $L$ via the Steiner quadruple systems $S(\mathcal C,v)$
associated to the codewords $v$ of $\mathcal C$. In this case, we
can take ${\mathcal C}/L$ as the vertex set of a quotient graph
$H_L({\mathcal C})$ of $M({\mathcal C})$ by setting an edge between
two classes $u+L$ and $v+L$ of ${\mathcal C}/L$ if and only if $uv$
is an edge of $M({\mathcal C})$.

\begin{propo}
Every extended 1-per\-fect code $\mathcal C$ is foldable over any
linear $L\subseteq K$. \qfd
\end{propo}

Recall that a covering graph map is a graph map $\phi:G\rightarrow G'$ for which
there is a nonnegative integer $s$ such that the inverse image
$\phi^{-1}$ of each vertex and of each edge of $G'$ has cardinality
$s$.

\begin{coro}
If $\mathcal C$ is foldable over a linear subspace $L$ of $K$, then
the natural projection ${\mathcal C}\rightarrow {\mathcal C}/L$ is
extendible to a covering graph map $\phi_L:M({\mathcal
C})\rightarrow H_L({\mathcal C})$. Moreover, if $\mathcal C$ is
foldable over $K$, then it is also foldable over $L$. \qfd
\end{coro}

In the setting of Corollary 4, given an edge $e=(v+L)(w+L)$ of
$H_L({\mathcal C})$, its {\it multiplicity} is the cardinality of
the set of labeling quadruples of edges in the SQS-subset
$\phi_L^{-1}(e)$. Note that the sum of the multiplicities of the
edges incident to any fixed vertex $v$ of $H_K({\mathcal C})$ must
equal 140, being this the cardinality of the SQS induced by
$\mathcal C$ at $v$. Of these 140 quadruples, 28 will be treated in
Theorem 5, item 1, and Section 6 in relation to the loops in each
$H_K({\mathcal C})$, where $\mathcal C$ is and SP-code with
$9\ge\kappa\ge 5$. The remaining 112 edges will be treated in
Theorem 5, item 2, and Section 7 as seven bunches of 16
quadruples formed as products of classes from partitions
$\{C_0,\ldots,C_7\}$ and $\{D_8,\ldots,D_{15}\}$.

\section{Participating 
STS(15)-types}

We denote the type of a Steiner triple
system of length 15, or STS(15)-type, by its associated integer
$t=1,\ldots,80$ in the common-order lists of the 80 existing
STS(15)-types in \cite{LeVan,MPR,WCC}. Pasch configurations, nicknamed as {\it fragments} in
\cite{LeVan}, are present in varying proportions in these STS(15)-types. An
algorithmic approach, (see also
\cite{Dej1,Dej2}), uses these fragments in order to determine the STS-graph
invariant $H_K(C)$. This yields the STS(15)-types
${\cal S}[v]$ associated with each one of the $2^{11}=2048$
codewords $v\in C$.

In the case of extended 1-perfect codes $\mathcal C$ of length 16,
first we obtain the STS-graphs modulo $K$ of the 16 punctured codes
${\mathcal C}_i, (i=1,\ldots,16),$ of length 15 that can be obtained
from each such code $\mathcal C$. Each such punctured code yields a
collection of 2048 16-tuples, determining, as mentioned, an integer in $[1,80]$
representing its STS(15)-type.

The integers in $[1,80]$ representing the punctured codes of
SP-codes with $9\geq\kappa\geq 5$ are: 2,3,4,5,6,7,8,13,14,16. According to \cite{LeVan}, they
have numbers of fragments accompanied by corresponding 15-tuples of
numbers of fragments containing each a specific coordinate index,
but given in nondecreasing order, as follows, (where the first line,
cited just for reference, stands for the linear code):
$$\begin{array}{rr}
_{1:} & _{105(42,42,42,42,42,42,42,42,42,42,42,42,42,42,42);} \\
^{2:}_{3:} & ^{73(42,30,30,30,30,30,30,30,30,26,26,26,26,26,26);}_{57(26,26,26,24,24,24,24,24,24,24,24,18,18,18,18);} \\
^{4:}_{5:} & ^{49(30,26,22,20,20,20,20,18,18,18,18,18,18,14,14);}_{49(26,26,20,20,20,20,18,18,18,18,18,18,18,18,18);} \\
^{6:}_{7:} & ^{37(22,22,22,14,14,14,14,14,14,12,12,12,12,12,12);}_{33(18,18,18,12,12,12,12,12,12,12,12,12,12,12,12);} \\
^{\,\,8:}_{13:} & ^{37(18,18,18,15,15,15,15,14,14,14,14,14,14,14,10);}_{33(20,16,16,14,14,12,12,12,12,12,12,12,12,12,10);} \\
^{14:}_{16:}
&^{37(24,16,16,16,15,15,15,15,14,14,14,12,12,12,12);}_{49(21,21,21,21,21,21,21,21,18,18,18,18,18,18,18).}
\end{array}$$
A list of the 16-tuples representing the vertices of the SQS-graphs
for the treated SP-codes can be found in

\noindent http://home.coqui.net/dejterij/tuples.txt.

\section{Which SP-codes $\mathcal C$ are homogeneous?}

An extended 1-perfect code $\mathcal C$ is {\it
SQS-homogeneous} if and only if the SQSs determined  by its
codewords are all equivalent: in case $\mathcal C$ is
of length 16, the 16 punctured codes of $\mathcal C$ have the same
distribution of STS(15)-types.

An SQS-homogeneous 1-perfect code of length 16 is {\it
STS-homogeneous} if and only if each one of its punctured codes is
homogeneous via a {\it common} STS(15). Among the SP-codes, we calculated that those
behaving in this fashion are, for each $\kappa$ such that
$9\geq\kappa\geq 5$, the following ones:
$$\begin{array}{ll}
^{\kappa=9:}_{\kappa=8:} & ^{(007,2), (008,2);}_{(114,5), (115,3), (963,g);} \\
^{\kappa=7:}_{\kappa=6:} & ^{(002,2), (003:2), (004,2);}_{(064,4), (917,8), (918,8);} \\
^{\kappa=5:} & ^{(708,d);}
\end{array}$$
where each code denomination is accompanied (between parentheses) by
the STS(15)-type, from 1 to 80, common as induced STS(15) to the
corresponding 2048 codewords, and the types 13, 14 and 16 are
respectively represented by the letters $c$, $d$ and $g$.

We also calculated that the remaining SQS-homogeneous SP-codes such that $9\ge\kappa\ge 5$ and having 16 STS-homogeneous punctured codes are:

$$\begin{array}{llll}
^{\kappa=8:} & ^{(112,3333777733337777),}_{(116,5555222233553355),}
 & ^{(113,5555555577557755),}_{(117,5555555522222222),} \\
_{\kappa=7:} & ^{(118,3333333333332222);}_{(101,4444222244444444),}
 & _{(102,4444444433223322),} \\
 & ^{(103,2244224422442244),}_{(105,4444444422332225),}
 & ^{(104,4422442222552255),}_{(959,8855885588558855),} \\
_{\kappa=6:} & ^{(960,8888888833gg33gg);}_{(063,4444444477557755),}
 & _{(065,5555555544444444),} \\
 & ^{(066,4444222244444444),}_{(068,4444222255335533),}
 & ^{(067,5555222244444444),}_{(919,8888444488884444),} \\
 & ^{(921,4488448844884488),}_{(923,55554444dd55dd55),}
 & ^{(922,8855885544884488),}_{(924,44444444dd55dd55),} \\
 & ^{(925,33333333dddddddd)}_{(931,8888222233ee33ee);}
 & ^{(930,gggg2222eeeeeeee),} \\
^{\kappa=5} & ^{(701,4444444477557755),}_{(706,dddd7777dddddddd),}
 & ^{(702,4444444477557755),}_{(709,88888888dddddddd),} \\
 & ^{(710,88888888dddddddd),}_{(706,dddd7777dddddddd),}
 & ^{(714,eeeeeeeeeeee5555),}_{(716,88885555eeee5555),} \\
 & ^{(717,88884444eeee5555),}_{(720,55558888ddeeddee),}
 & ^{(719,88884444ddeeddee),}_{(721,55554444ddee4444),} \\
 & ^{(722,55554444ddee4444),}_{(726,3333dddddddd3333);}
 & ^{(725,eeee5555eeee5555),} \\
\end{array}$$
where each code denomination is accompanied between parenthesis by
the common 16-tuple of STS(15)-types formed from the 16 punctured
codes.

We also checked that there are still some SQS-homogeneous SP-codes with $9\ge\kappa\ge 5$ whose punctured codes
{\it are not} STS-ho\-mo\-ge\-neous:
$$\begin{array}{llll}
^{\kappa=6:} & ^{(914,333388888888gggg),}_{(916,333355558888gggg),}
 & ^{(915,333388888888gggg),}_{(926,223333ddddddddgg),} \\
 & ^{(927,223333ddddddddgg),}_{(929,224444445588dddd);}
 & ^{(928,224444445588dddd),} \\
^{\kappa=5:} & ^{(029,3344445555666677),}_{(705,44448888ddddeeee),}
 & ^{(704,44444444ddddeeee),}_{(707,55777777dddddddd),} \\
 & ^{(711,44448888dddddddd),}_{(718,222344448888gggg),}
 & ^{(713,44555588eeeeeeee),}_{(723,22333344445555gg);} \\
\end{array}$$
where STS(15)-type denomination numbers are given in non-decreasing
order, but their actual order differs coordinate by coordinate. For
example, SP-code ${\mathcal C}=914$, which is SQS-homogeneous, has
$H_K({\mathcal C})$ holding 16 vertices yielding the 16-tuple
$888888883g3g3g3g$ and 16 vertices yielding the 16-tuple
$88888888g3g3g3g3\neq 888888883g3g3g3g$.

\section{What do the edges of $H_K({\mathcal C})$ stand for?}

In what follows we present Theorem 5, accounting for the structure of
the SQS-subsets $\phi_K^{-1}(e)$ represented by the edges $e$ of
$H_K({\mathcal C})$, for the 361 SP-codes $\mathcal C$ with
$9\geq\kappa\geq 5$. (Recall that $|\phi_K^{-1}(e)|$ is the
multiplicity of $e$). To express the coordinate indices in codewords
of $\mathcal C$, we use hexadecimal notation: these indices
constitute the set $[0,f]=\{0,1,\ldots,9,a,b,\ldots,f\}$.

Let $0<s\in\Z$. If $Y$ is a set of quadruples of $[0,s]$, let the
$s$-{\it supplement} of $Y$ be the set of quadruples
$\{x_1,x_2,x_3,x_4\}\in[0,s]$ such that
$\{s-x_1,s-x_2,s-x_3,s-x_4\}\in Y.$ 

Let $S=\{s_1,\ldots,s_t\}$ be a partition of number $7$
into positive integers $s_i$ such that $s=s_1+\ldots +s_t$ and
$s_1\leq\ldots\leq s_t$. Let $Y$ be a $t$-set of quadruples of
$[0,7]$. We define a {\it descending}, (resp. an {\it ascending}),
$S$-{\it partition} ${\mathcal P}_Y^\downarrow$, (resp. ${\mathcal
P}_Y^\uparrow$), to be a partition $\{Y_1,\ldots,Y_t\}$ of $Y$
such that $|Y_i|=s_i$, for each $i\in[1,t]$, and if $(w_i,w_j)\in Y_i\times Y_j$, where $i,j\in[1,t]$ and $i<j$, then $w_i>w_j$, (resp.
$w_i<w_j$), lexicographically.

If $S$ has $s_1$, (resp. $s_t)$, equal to $\min\{2^{\kappa-5}-1,7\}$ and has every
other $s_i$ equal to $\min\{2^{\kappa-5},7\}$, then we say that
${\mathcal P}_Y^\downarrow$, (resp. ${\mathcal P}_Y^\uparrow$), is a
{\it descending}, (resp. {\it ascending}), $(\kappa-5)$-{\it
partition}.

Associated to the Fano plane on vertex set $[1,7]$ and line set
$\{123,145,$ $167,247, 256,346,357\}$, we have the following three
sets of quadruples:
$$\begin{array}{ll}
X= & \{0123, 0145, 0167, 0247, 0256, 0346, 0357\} \\
Y= & \{4567,2367,2345,1356,1347,1257,1246\} \\
Z= &
\{       cdef,abef,89ef,\,8bdf,\,9adf,\,9bcf,\,8acf, \\
 & \,\,\,89ab,89cd,abcd,\,9ace,\,8bce,\,8ace,\,9bce\},
\end{array}$$
where 
$Z=f$-supplement of $X\cup Y$.

\begin{thm}
Let $\mathcal C$ be an SP-code of length 16 with $9\ge\kappa\ge 5$. Then:
\begin{enumerate}\item
each vertex $v$ of $H_K({\mathcal C})$ has a loop $\ell_v$ of
multiplicity $|\phi_K^{-1}(\ell_v)|$, with its SQS-subset
$\phi_K^{-1}(\ell_v)$ formed as the union of:

\noindent{\bf(a)} $Z$; {\bf(b)}
the last set $Y_t$ in the descending $(\kappa-5)$-partition
${\mathcal P}_Y^\downarrow$, where $t=2^{\max\{0,{7-\kappa}\}}$;
{\bf(c)} $X$, if $\kappa\geq 8$; {\bf(d)} a specific product
$C_{i_0}\times D_{8+j_0}$ of partition classes $C_{i_0}$ and $D_{8+j_0}$,
if $\kappa=9$;
\item[\bf 2.]
each link $e$ of $H_K({\mathcal C})$ has $\phi_K^{-1}(e)$ formed as the union of:

\noindent {\bf(a)}
a union of lexicographically ordered quarters (LOQs) of products $C_i\times
D_{8+j}$ ($\neq C_{i_0}\times D_{8+j_0}$ in item 1 above, if $\kappa=9$), namely:
{\bf(i)}
the eight LOQs in two such products, if $\kappa=9$;
{\bf(ii)}
the four LOQs in one such product, if $\kappa=8$;
{\bf(iii)} between one and three
LOQs in one such product, if $\kappa\leq 7$;
{\bf(b)}
at most either one set $\neq Y_t$ in the
descending $(\kappa-5)$-partition ${\mathcal P}_Y^\downarrow$ or one
set in the ascending $(\kappa-5)$-partition ${\mathcal
P}_X^\uparrow$, if $\kappa\leq 7$.
\end{enumerate}
\end{thm}

\proof The properties of the loops, resp. links, of $\mathcal C$,
establishing the statement of the theorem for the five treated kernel dimensions $\kappa$
are considered in Section 6, resp. 7, with data details for the 361 cases considered present in http://home.coqui.net/dejterij/details.tex.
In Section 6, the diagonal and nondiagonal elements in the tables closing the three final subsections, (that is, for $\kappa\leq 7$), allow to establish items 1(b) and 2(b), respectively, in the statement of the theorem.

\qfd

\section{Vertices and loops of $H_K({\mathcal C})$}

\subsection{Case $\kappa=9$}

For each SP-code $\mathcal C$ with $\kappa=9$, namely for ${\mathcal
C}=007$ and 008, there is a subspace $L$ of index 2 in
$K=Ker({\mathcal C})$ such that the vertices of $H_L({\mathcal C})$
are given by eight classes mod $L$ that we denote $k=0,\ldots,7$, leading to four classes mod $K$ formed by the union of classes
$2j$ and $2j+1$ mod $L$, for $j=0,1,2,3$. This graph $H_L({\mathcal C})$ has a loop of multiplicity 28 at each vertex of $H_L({\mathcal C})$ represented by
$X\cup Y\cup Z$. This loop together with an edge of
multiplicity 16 obtained from a product as in the table of
Subsection 7.1 below, for each vertex of $H_L({\mathcal C})$,
project onto a loop of multiplicity $28+16=44$ in $H_K({\mathcal
C})$.

\subsection{Case $\kappa=8$}

The vertices of each one of the eight existing $H_K({\mathcal C})$
here, namely for
$${\mathcal C}=005, 006, 112, 113, 114, 115, 116, 963,$$ are given by 8 classes
mod $K$ that we denote $k=0,\ldots,7$. Each such class has a loop
of multiplicity 28, represented by $X\cup Y\cup Z$.

\subsection{Case $\kappa=7$}

The vertices of each one of the 18 existing $H_K({\mathcal C})$
here, namely for
$${\mathcal C}=002,\ldots,004,101,\ldots,111,959,\ldots,962,$$ are given by 16
classes mod $K$ that we denote $k_n$, with $k=0,\ldots,7$ and $n=0,1$. Each
$H_K({\mathcal C})$ presents a loop of multiplicity 21 represented
by $X'=Y\cup Z$ and an edge of multiplicity 7 represented by $X$.
The following contributive table for SQS-subsets $\phi_K^{-1}(\ell)$
of loops and links $\ell$ of $H_K({\mathcal C})$ holds, with
multiplicities indicated between parenthesis:
$$\begin{array}{c|cc|}
 &    k_0 & k_1 \\ \hline
 k_0 & X'(21) & X(7) \\
 k_1 & X(7) &  X'(21)
\end{array}$$

\subsection{Case $\kappa=6$}

The vertices of each one of the 86 existing $H_K({\mathcal C})$
here, namely for
$${\mathcal C}=063,064,065,\ldots,098,099, 100, 911,912,913,\ldots,956,957,958,$$
are given by 32 classes mod $K$ that we denote $k_n$, where $k=0,\ldots,7$ and $n=0,1,2,3$.
Let
$$\begin{array}{ll}
A=\{0123,0145,0167\}, & B=\{0247,0256,0346,0357\}, \\
A'=[0,7]\mbox{-supplement of }A, & B'=[0,7]\mbox{-supplement of }B,
\end{array}$$
and $Z'=Z\cup A'$. The following contributive table for SQS-subsets
$\phi_K^{-1}(\ell)$ of loops and links $\ell$ of $H_K({\mathcal C})$
holds, with multiplicities indicated between parenthesis:
$$\begin{array}{c|cccc|} & k_0  & k_1  & k_2  & k_3  \\ \hline
 k_0 & Z'(17) & B'(4) &  B(4) & A(3) \\
 k_1 & B'(4) & Z'(17) & A(3) & B(4) \\
 k_2 & B(4) & A(3) & Z'(17) & B'(4) \\
 k_3 & A(3) & B(4) & B'(4) & Z'(17)
 \end{array}$$

\subsection{Case $\kappa=5$}

The vertices of each one of the 244 existing \ $H_K({\mathcal C})$
here, namely for
$${\mathcal C}=029,\ldots,060,061, 062, 701,702,703,704,\ldots,906,907,908,909,910,$$
are given by 64 classes mod $K$l that we denote $k_n$, where
$k,n\in\{0,\ldots,7\}$. Let $A_0=\{0123\}$,  $A_1=\{0145,0167\}$,
$B_0=\{0247,0256\}$,  $B_1=\{0346,$ $0357\}$; $A_i'=[0,7]$-supplement
of $A_i$,   $B_i'=[0,7]$-supplement of $B_i$, for $i=0,1$, and
$Z_0=Z\cup A_0'$. The following contributive table for SQS-subsets
$\phi_K^{-1}(\ell)$ of loops and links $\ell$ of $H_K({\mathcal C})$
holds, with multiplicities indicated between parenthesis and $f=15$:
$$\begin{array}{c|cccccccc}
&   k_0  & k_1 & k_2 & k_3 & k_4 & k_5 & k_6 & k_7 \\ \hline
k_0 & Z_0(f) & A_1'(2) & B_0'(2) & B_1'(2) & B_1(2) & B_0(2) & A_1(2) & A_0(1) \\
k_1 & A_1'(2) &  Z_0(f) & B_1'(2) & B_0'(2) &  B_0(2) & B_1(2) & A_0(1) & A_1(2) \\
k_2 & B_0'(2) & B_1'(2) & Z_0(f) & A_1'(2) & A_1(2) &  A_0(1) & B_1(2) & B_0(2) \\
k_3 & B_1'(2) & B_0'(2) & A_1'(2) & Z_0(f) & A_0(1) &  A_1(2) & B_0(2) &  B_1(2) \\
k_4  & B_1(2) & B_0(2) & A_1(2) & A_0(1) & Z_0(f) & A_1'(2) & B_0'(2) & B_1'(2) \\
k_5 & B_0(2) & B_1(2) & A_0(1) & A_1(2) & A_1'(2) & Z_0(f) & B_1'(2) & B_0'(2) \\
k_6 & A_1(2) & A_0(1) & B_1(2) & B_0(2) & B_0'(2) & B_1'(2) & Z_0(f) & A_1'(2) \\
k_7 & A_0(1) & A_1(2) & B_0(2) & B_1(2) & B_1'(2) & B_0'(2) & A_1'(2) & Z_0(f)
\end{array}$$

\section{Links of $H_K({\mathcal C})$}

In this section, we specify the form of the products claimed in
Theorem 5. The actual denomination numbers in
$\{0,\ldots,6,8,9,10\}$ for the partitions $\{C_0,\ldots,C_7\}$ and
$\{D_8,\ldots,D_f\}$ of Theorem 1, which are used in those
pro\-ducts, for each SP-code $\mathcal C$ with $9\ge\kappa\ge 5$, are
shown explicitly in http://home.coqui.net/dejterij/details.pdf.

\subsection{Case $\kappa=9$}

Consider the following partitions of length 7:
$$\begin{array}{ccccccc}
1_a=1_3^5, & 2_a=2_3^7, & 3_a=3_2^7, & 4_a=4_5^7, & 5_a=5_4^6, & 6_a=6_7^4, & 7_a=7_6^5, \\
1_b=1_3^6, & 2_b=2_3^6, & 3_b=3_2^6, & 4_b=4_5^6, & 5_b=5_4^7, &
6_b=6_7^5, & 7_b=7_6^4,
\end{array}$$
where two notations for each partitions are used. The first
notation, $\alpha_p$, where $\alpha=0,\ldots,7$ and $p$ is a letter,
is a shorthand used in the tables below. The second notation,
$\alpha_\ell^m$, represents the partition with a lexicographically ordered
form $\{0\alpha,x\ell,ym,zw\}$, as in the subsequent example, where the symbol $\alpha_p\beta_q$ ($=1_a1_a$), with
$\alpha,\beta\in[0,7]$ and $p,q\in\{a,b\}$, represents the
product $\alpha_p\times(\beta+8)_q$ of the two partitions $\alpha_p$
and $(\beta+8)_q$, in each case:
\begin{eqnarray}
1_a1_a=1_3^51_3^5=\{01,23,45,67\}\times\{89,ab,cd,ef\}=\{0189, \ldots,
67ef\}. 
\end{eqnarray}
In case ${\mathcal C}=007$, we have the following contributive table
for SQS-subsets $\phi_L^{-1}(e)$ of links $e$ of $H_L({\mathcal
C})$, where $L$ is as in Subsection 6.1:
$$\begin{array}{l|cccccccc}
k\setminus\ell & 0 & 1 & 2 & 3 & 4 & 5 & 6 & 7 \\ \hline
0  & ....   & 1_a2_a & 3_b1_a & 2_a3_b & 5_a5_a & 4_a7_b & 6_b6_b & 7_b4_a \\
1  & 1_a2_a & ....   & 2_a3_b & 3_b1_a & 4_a6_b & 5_a4_a & 7_b5_a & 6_b7_b \\
2  & 3_b1_a & 2_a3_b & ....   & 1_a2_a & 7_b4_a & 6_b6_b & 4_a7_b & 5_a5_a \\
3  & 2_a3_b & 3_b1_a & 1_a2_a & ....   & 6_b7_b & 7_b5_a & 5_a4_a & 4_a6_b \\
4  & 5_a5_a & 4_a6_b & 7_b4_a & 6_b7_b & ....   & 1_a3_b & 2_a2_a & 3_b1_a \\
5  & 4_a7_b & 5_a4_a & 6_b6_b & 7_b5_a & 1_a3_b & ....   & 3_b1_a & 2_a2_a \\
6  & 6_b6_b & 7_b5_a & 4_a7_b & 5_a4_a & 2_a2_a & 3_b1_a & ....   & 1_a3_b \\
7  & 7_b4_a & 6_b7_b & 5_a5_a & 4_a6_b & 3_b1_a & 2_a2_a & 1_a3_b &
....
\end{array}$$
In this table, the 16 quadruples corres\-ponding to each
sub-diagonal entry form the product $C_i\times C_j$ contributing to
the SQS-subset $\phi^{-1}(\ell)$ of a corresponding loop $\ell$ of
$H_K({\mathcal C})$ as in item 1 of Theorem 5. The SQS-subsets
$\phi_K^{-1}(e)$ for links $e$ of $H_K({\mathcal C})$ are obtained by
considering that the vertices of $H_K({\mathcal C})$ are unions of
the classes $2j$ and $2j+1$ mod $L$, for $j=0,1,2,3$.

A similar disposition for the case 008 is shown in tabulated format
in http://home.coqui.net/dejterij/$xyz$PAT.txt, where $xyz=008$. By
replacing $xyz$ by any other 3-string of an SP-code with
$9\ge\kappa\ge 5$, a corresponding file may be downloaded.

\subsection{Case $\kappa=8$}

We deal here with 8 classes mod $K$, (instead of 8 classes mod $L$,
as above). For ${\mathcal C}=005$, we have the following
contributive table for SQS-subsets $\phi_K^{-1}(e)$ of links $e$ of
$H_K({\mathcal C})$, (otherwise, we refer to the last comment in
Subsection 7.1):

$$\begin{array}{l|cccccccc}
k\setminus\ell &    0 &   1 & 2  &  3 &  4 & 5  & 6  & 7 \\ \hline
0 & .... & 1_a1_a & 3_b3_b & 2_a2_a & 5_a5_a & 4_a4_a & 6_a6_a & 7_a7_a \\
1 & 1_a1_a & .... & 2_a2_a & 3_b3_b & 4_a4_a & 5_a5_a & 7_a7_a & 6_a6_a \\
2 & 3_b3_b & 2_a2_a & .... & 1_a1_a & 7_b7_b & 6_b6_b & 4_b4_b & 5_b5_b \\
3 & 2_a2_a & 3_b3_b & 1_a1_a & .... & 6_b6_b & 7_b7_b & 5_b5_b & 4_b4_b \\
4 & 5_a5_a & 4_a4_a & 7_b7_b & 6_b6_b & .... & 1_a1_a & 2_b2_b & 3_a3_a \\
5 & 4_a4_a & 5_a5_a & 6_b6_b & 7_b7_b & 1_a1_a & .... & 3_a3_a & 2_b2_b \\
6 & 7_b7_b & 6_b6_b & 5_a5_a & 4_a4_a & 3_b3_b & 2_a2_a & .... & 1_a1_a \\
7 & 6_b6_b & 7_b7_b & 4_a4_a & 5_a5_a & 2_a2_a & 3_b3_b & 1_a1_a &
....
\end{array}$$

\subsection{Case $\kappa=7$}

In addition to the partitions mentioned in the subsections above, we
need the following ones:
$$\begin{array}{ccccccc}
1_c=1_3^7, & 2_c=2_3^5, & 3_c=3_2^5, & 4_c=4_6^5, & 4_d=4_6^7, & 4_e=4_7^6, & 5_c=5_7^4, \\
5_d=5_7^6, & 5_e=5_6^7, & 6_c=6_4^7, & 6_d=6_4^5, & 6_e=6_5^4, & 7_c=7_5^6, & 7_d=7_5^4, \\
7_e=7_4^5. &            &            &            &            & &
\end{array}$$
Codes 002, 003, 004, 102, 103, 104, 106, 107, 109, 959, 960, 961,
962, (resp. 101, 105), [resp. 108, 110, 111], use partitions of the
form $\alpha_a,\alpha_b$, (resp. $\alpha_c,\alpha_d$), [resp.
$\alpha_c,\alpha_e$], where $\alpha=0,\ldots,7$.

For example, in the case 002, we have the following contributive
table for SQS-subsets $\phi_L^{-1}(e)$ of links $e$ of
$H_K({\mathcal C})$, where $m=0,1$:
$$\begin{array}{l|cccccccc}
k_n\setminus\ell_m & 0_m & 1_m & 2_m & 3_m & 4_m & 5_m & 6_m & 7_m \\
\hline
0_0 & .... & 1_a1_a & 2_b7_a & 3_a6_a & 5_a3_a & 4_a2_b & 6_a4_a & 7_a5_a \\
0_1 & .... & 1_a1_a & 3_a7_a & 2_b6_a & 5_a3_a & 4_a2_b & 7_a4_a & 6_a5_a \\
1_0 & 1_a1_a & .... & 3_a6_a & 2_b7_a & 4_a2_b & 5_a3_a & 7_a5_a & 6_a4_a \\
1_1 & 1_a1_a & .... & 2_b6_a & 3_a7_a & 4_a2_b & 5_a3_a & 6_a5_a & 7_a4_a \\
2_0 & 3_b6_b & 2_a7_b & .... & 1_a1_a & 7_b4_b & 6_b5_b & 4_b3_b & 5_b2_a \\
2_1 & 2_a6_b & 3_b7_b & .... & 1_a1_a & 6_b4_b & 7_b5_b & 4_b3_b & 5_b2_a \\
3_0 & 2_a7_b & 3_b6_b & 1_a1_a & .... & 6_b5_b & 7_b4_b & 5_b2_a & 4_b3_b \\
3_1 & 3_b7_b & 2_a6_b & 1_a1_a & .... & 7_b5_b & 6_b4_b & 5_b2_a & 4_b3_b \\
4_0 & 5_a2_a & 4_a3_b & 6_a4_b & 7_a5_b & .... & 1_a1_a & 2_b7_b & 3_a6_b \\
4_1 & 5_a2_a & 4_a3_b & 7_a4_b & 6_a5_b & .... & 1_a1_a & 3_a7_b & 2_b6_b \\
5_0 & 4_a3_b & 5_a2_a & 7_a5_b & 6_a4_b & 1_a1_a & .... & 3_a6_b & 2_b7_b \\
5_1 & 4_a3_b & 5_a2_a & 6_a5_b & 7_a4_b & 1_a1_a & .... & 2_b6_b & 3_a7_b \\
6_0 & 7_b4_a & 6_b5_a & 4_b2_b & 5_b3_a & 3_b6_a & 2_a7_a & .... & 1_a1_a \\
6_1 & 6_b4_a & 7_b5_a & 4_b2_b & 5_b3_a & 2_a6_a & 3_b7_a & .... & 1_a1_a \\
7_0 & 6_b5_a & 7_b4_a & 5_b3_a & 4_b2_b & 2_a7_a & 3_b6_a & 1_a1_a & .... \\
7_1 & 7_b5_a & 6_b4_a & 5_b3_a & 4_b2_b & 3_b7_a & 2_a6_a & 1_a1_a &
....
\end{array}$$

The LOQs in which the
products $\alpha_p\beta_q$ in the table above divide are the
destinations of the classes $k_n$ in $M({\mathcal C})$ that yield
the contributions to the SQS-subsets $\phi_K^{-1}(e)$ of the edges
$e$ of $H_K({\mathcal C})$. A similar second table can be set with a
symbol $\epsilon_1\epsilon_2\epsilon_3\epsilon_4$ in each
non-diagonal entry, each $\epsilon_i$ representing a LOQ of an
$\alpha_p\beta_q$. In fact, for $k_n$ with $n=0$, ($n=1$), we have:
$\epsilon_1\epsilon_2\epsilon_3\epsilon_4=1000$,
($\epsilon_1\epsilon_2\epsilon_3\epsilon_4=0111$), where
$k\in[0,7]$. For example, $1_3^51_3^5$ in position
$(k_n,\ell_m)=(0_0,1_m)$ in the table above has
$\epsilon_1\epsilon_2\epsilon_3\epsilon_4=1000$ in this second
table, meaning that $0_0$ assigns
$$\begin{array}{c}
\{018a, 019b, 01cd, 01ef\}\mbox{ to }\ell_m=\ell_{\epsilon_1}=1_0; \\
\{238a, 239b, 23cd, 23ef\}\mbox{ to }\ell_m=\ell_{\epsilon_2}=0_0; \\
\{458a, 459b, 45cd, 45ef\}\mbox{ to }\ell_m=\ell_{\epsilon_3}=0_0; \\
\{678a, 679b, 67cd, 67ef\}\mbox{ to }\ell_m=\ell_{\epsilon_4}=0_0.
\end{array}$$
A listing showing the combination of the file $xyz$PAT.txt mentioned
in Subsection 7.1 and the $\epsilon_1\epsilon_2\epsilon_3\epsilon_4$
above can be retrieved from

\noindent http://home.coqui.net/dejterij/$xyz$TEST.txt, where $xyz$ is 002 or
the value of $xyz$ corresponding to any  SP-code with $\kappa=7,6,5$.

\subsection{Case $\kappa=6$}

Consider the partitions of length 7 given above together with:
$$\begin{array}{cccccc}
1_d=1_5^{4}, & 1_e=1_6^{7}, & 1_f=1_4^{5}, & 2_d=2_5^{7}, & 2_e=2_4^{6}, & 2_f=2_6^{4}, \\
3_d=3_4^{7}, & 3_e=3_7^{4}, & 3_f=3_5^{6}, & 4_f=4_3^{6}, & 4_g=4_2^{7}, & 4_h=4_5^{3}, \\
5_f=5_2^{6}, & 5_g=5_3^{7}, & 6_f=6_2^{5}, & 6_g=6_7^{3}, & 7_f=7_6^{3}, & 7_g=7_3^{5}. \\
\end{array}$$
For each SP-code $\mathcal C$ with $\kappa=6$, we can assign a
product $\alpha_p\beta_q$ to each class $k_n=0_0,\ldots,7_3$ in
eight tables, each with rows headed by $k_n$, where $k\in[0,7]$ is
fixed and $n$ varies in $[0,3]$, and with columns headed by all
values of $k_m$, where $m\in[0,3]$ is fixed and $k$ varies in
$[0,7]$. Each of these eight tables expresses the needed products
$\alpha_p\beta_q$. We exemplify the values of $p$ for the case
${\mathcal C}=066$ in a table with each entry $(k,n)\in
[0,7]\times[0,3]$ containing a literal 7-tuple
$(p_1,p_2,p_3,p_4,p_5,p_6,p_7)$ which is the 7-tuple of subindexes
in a corresponding expression
$(1_{p_1},2_{p_2},3_{p_3},4_{p_4},5_{p_5},6_{p_6},7_{p_7})$:
$$\begin{array}{l|cccc}
_{k\setminus n} &   _0     &    _1    &   _2     &    _3    \\
\hline
^0_1 & ^{abacdcd}_{abadcdc} & ^{abacdcd}_{abadcdc} & ^{abadcdc}_{abacdcd} & ^{abadcdc}_{abacdcd} \\
^2_3 & ^{aabcdcd}_{aabdcdc} & ^{aabdcdc}_{aabcdcd} & ^{aabcdcd}_{aabdcdc} & ^{aabdcdc}_{aabcdcd} \\
^4_5 & ^{bacaaaa}_{bac\,bbbb} & ^{b\,ac\,bbbb}_{bacaaaa} & ^{bacaaaa}_{b\,ac\,bbbb} & ^{b\,ac\,bbbb}_{bacaaaa} \\
^6_7 & ^{b\,ca\,bbbb}_{bcaaaaa} & ^{b\,ca\,bbbb}_{bcaaaaa} & ^{bcaaaaa}_{b\,ca\,bbbb} & ^{bcaaaaa}_{b\,ca\,bbbb} \\
\end{array}$$
The components $\beta_q$ of products $\alpha_p\beta_q$ here are
constant-partition 4-tuples. The information for case ${\mathcal
C}=066$ can be condensed as follows:

$$\begin{array}{l|cccccccc}
k_n\setminus\ell_m & 0_m & 1_m & 2_m & 3_m & 4_m & 5_m & 6_m & 7_m \\
\hline
0_n & ....         & 1_{p_1}1_a & 3_{p_3}3_a & 2_{p_2}2_b & 4_{p_4}4_a & 6_{p_6}7_a & 5_{p_5}6_a & 7_{p_7}5_a \\
1_n & 1_{p_1}1_a & ....         & 2_{p_2}2_b & 3_{p_3}3_a & 5_{p_5}5_a & 7_{p_7}6_a & 4_{p_4}7_a & 6_{p_6}4_a \\
2_n & 3_{p_3}2_a & 2_{p_2}3_b & ....         & 1_{p_1}1_a & 6_{p_6}7_b & 4_{p_4}4_b & 7_{p_7}5_b & 5_{p_5}6_b \\
3_n & 2_{p_2}3_b & 3_{p_3}2_a & 1_{p_1}1_a & ....         & 7_{p_7}6_b & 5_{p_5}5_b & 6_{p_6}4_b & 4_{p_4}7_b \\
4_n & 5_{p_5}4_a & 4_{p_4}5_a & 7_{p_7}6_a & 6_{p_6}7_a & ....         & 1_{p_1}2_b & 2_{p_2}3_a & 3_{p_3}1_a \\
5_n & 7_{p_7}6_a & 6_{p_6}7_a & 5_{p_5}4_a & 4_{p_4}5_a & 1_{p_1}3_a & ....         & 3_{p_3}1_a & 2_{p_2}2_b \\
6_n & 4_{p_4}7_b & 5_{p_5}6_b & 6_{p_6}5_b & 7_{p_7}4_b & 2_{p_2}2_a & 3_{p_3}1_a & ....         & 1_{p_1}3_b \\
7_n & 6_{p_6}5_b & 7_{p_7}4_b & 4_{p_4}7_b & 5_{p_5}6_b & 3_{p_3}1_a
& 2_{p_2}3_b & 1_{p_1}2_a & ....
\end{array}$$
where $n,m\in[0,3]$. An observation on LOQs
of the $\alpha_p\beta_q$-s si\-mi\-lar to the one in Subsection 7.3
holds here. In fact, an accompanying table for the resulting
4-tuples $\epsilon_1\epsilon_2\epsilon_3\epsilon_4$ of LOQs can be
composed by replacing the symbols $A$ and $B$ in the simplified
table on the left side below (accompanying the one above) by the
sub-tables on its right side:
$$\begin{array}{l|cccccccc||cc|cc||c|c}
_{k\setminus\ell} & _0 & _1 & _2 & _3 & _4 & _5 & _6 & _7 & & _A &
_{\epsilon_1\epsilon_2\epsilon_3\epsilon_4} &  & _B &
_{\epsilon_1\epsilon_2\epsilon_3\epsilon_4} \\ \hline
^0_1 & ^._A & ^A_. & ^B_B & ^B_B & ^B_B & ^B_B & ^B_B & ^B_B & & ^{k_0}_{k_1} & ^{3000}_{2111} & & ^{k_0}_{k_1} & ^{2100}_{3110} \\
^2_3 & ^B_B & ^B_B & ^._A & ^A_. & ^B_B & ^B_B & ^B_B & ^B_B & & ^{k_2}_{k_3} & ^{1222}_{0333} & & ^{k_2}_{k_3} & ^{0322}_{1332} \\
^4_5 & ^B_B & ^B_B & ^B_B & ^B_B & ^._A & ^A_. & ^B_B & ^B_B &  &     &      &  &     &      \\
^6_7 & ^B_B & ^B_B & ^B_B & ^B_B & ^B_B & ^B_B & ^._A & ^A_. &  &     &      &  &     &      \\
\end{array}$$

\subsection{Case $\kappa=5$}

In addition to the partitions given above, consider:
$$\begin{array}{ccccccc}
1_g=1_4^7, & 1_h=1_4^6, & 1_i=1_7^6, & 1_j=1_5^7, & 1_k=1_6^5, & 2_g=2_7^6, & 2_h=2_5^6,  \\
2_i=2_4^7, & 2_j=2_7^5, & 2_k=2_4^5, & 3_g=3_7^6, & 3_h=3_7^5, & 3_i=3_4^6, & 3_j=3_6^7,  \\
4_i=4_2^6, & 4_j=4_3^5, & 5_h=5_6^4, & 5_i=5_4^3, & 5_j=5_2^4, & 5_k=5_6^3, & 6_h=6_3^5,  \\ 
6_i=6_3^4, & 6_j=6_5^7, & 6_k=6_5^3, & 7_h=7_5^3, & 7_i=7_2^5, &
7_j=7_3^4, & 7_k=7_2^3.
\end{array}$$
For each SP-code $\mathcal C$ with $\kappa=5$, we can assign a
product $\alpha_p\beta_q$ to each class $k_n=0_0,\ldots,7_7$ in
eight tables, each with rows headed by $k_n$, where $k\in[0,7]$ is
fixed and $n$ varies in $[0,7]$, and with columns headed by all the
values of $k_m$, where $m\in[0,7]$ is fixed and $k$ varies in
$[0,7]$. Each of these eight tables expresses the needed products
$\alpha_p\beta_q$. A condensed form of these eight tables exists as
in Subsection 7.4 and can be obtained from the sources cited at the
end of Subsections 7.1 and 7.3.

An observation on LOQs
of the $\alpha_p\beta_q$-s si\-mi\-lar to those in Subsections 7.3-4
holds. An accompanying table for the resulting 4-tuples
$\epsilon_1\epsilon_2\epsilon_3\epsilon_4$ of LOQs can be composed
by replacing the symbols $A,B,C,D$ in the simplified table below, to
the left (accompanying the one above) by the sub-tables on its right
side, where the column headers $A,B,C,D$ stand for the corresponding
4-tuples $\epsilon_1\epsilon_2\epsilon_3\epsilon_4$ of LOQs:
$$\begin{array}{l|cccccccc||c||c|c|c|c|}
_{k\setminus\ell}
  & _0 & _1 & _2 & _3 & _4 & _5 & _6 & _7 & & _A   & _B   & _C  & _D   \\ \hline 
 ^0_{1} & ^{.}_{A} & ^{A}_{.} & ^{B}_{C} & ^{C}_{B} & ^{D}_{D} & ^{D}_{D} & ^{D}_{D} & ^{D}_{D} & ^{k_0}_{k_1} & ^{6100}_{7110} & ^{4300}_{5211} & ^{5200}_{4311} & ^{4210}_{5310} \\
 ^2_{3} & ^{B}_{C} & ^{C}_{B} & ^{.}_{A} & ^{A}_{.} & ^{D}_{D} & ^{D}_{D} & ^{D}_{D} & ^{D}_{D} & ^{k_2}_{k_3} & ^{4322}_{5332} & ^{6221}_{7330} & ^{7220}_{6331} & ^{6320}_{7321} \\
 ^4_{5} & ^{D}_{D} & ^{D}_{D} & ^{D}_{D} & ^{D}_{D} & ^{.}_{A} & ^{A}_{.} & ^{C}_{B} & ^{B}_{C} & ^{k_4}_{k_5} & ^{2544}_{3554} & ^{0744}_{1655} & ^{1644}_{0755} & ^{0654}_{1754} \\
 ^6_{7} & ^{D}_{D} & ^{D}_{D} & ^{D}_{D} & ^{D}_{D} & ^{C}_{B} & ^{B}_{C} & ^{.}_{A} & ^{A}_{.} & ^{k_6}_{k_7} & ^{0766}_{1776} & ^{2665}_{3774} & ^{3664}_{2775} & ^{2764}_{3765} \\
\end{array}$$

\end{document}